\documentstyle[12pt,amssymb]{amsart}
\newcommand{\PP}{{\mathbb P}}
\newcommand{\R}{{\mathbb R}}
\newcommand{\C}{{\mathbb C}}

\newcommand{\N}{{\mathbb N}}
\newcommand{\T}{{\mathbb T}}
\newcommand{\Pnd}{{\mathbb P}_{n,d}}
\newcommand{\HE}{{\cal HE}_{n,d}}
\newcommand{\HH}{{\cal H}_{n,d}}

\renewcommand{\phi}{\varphi}

\begin{document}

\centerline{\underline{\bf Mikhail Zaidenberg}}

\bigskip

\centerline{\underline{\bf HYPERBOLIC SURFACES in $\PP^3$ : EXAMPLES}}

\bigskip

\centerline{\underline{\bf (after a joint work with B. Shiffman)}}

\bigskip

\centerline{\underline{\bf Carcassonne, 7-9.11.2003}}

\vskip 0.3in

{\footnotesize \bf\centerline{\underline{CONTENTS:}}

\bigskip

\begin{enumerate}
\item[1.] Generalities \hfill 1
\item[2.] Kobayashi Problems \hfill 3
\item[] Hyperbolicity of complements \hfill 3
\item[] Hyperbolicity of projective hypersurfaces \hfill 4
\item[] Methods employed \hfill 4
\item[] Examples of hyperbolic surfaces in $\PP^3$ \hfill 5
\item[3.] Symmetric powers of curves as hyperbolic hypersurfaces in
$\PP^3$ and $\PP^4$ \hfill
6
\item[4.] Deformation method and more hyperbolic surfaces in $\PP^3$ \hfill
8
\item[5.] Algebraic hyperbolicity \hfill 10
\item[6.] Bloch Conjecture \hfill 11
\item[] References \hfill 12
\item[] Appendix \hfill 14-25
\end{enumerate}}

\vskip 0.2in

\centerline{\underline{\bf 1. Generalities}}

\vskip 0.2in

{  \noindent {\bf DEFINITION}: On any complex space $X$ there
exists a unique \\ \smallskip {\it KOBAYASHI PSEUDOMETRIC}
$k_X\,:\,X\times X \to \R$ such that:

\medskip

(i) If $X=\Delta$ is the unit disc then $k_X$ is the Poincar\'e
metric;

\medskip

(ii) any holomorphic map $\varphi:\Delta \to X$ is a contraction:
$\varphi^*(k_X)\le k_\Delta$, and

\medskip

(iii) $k_X$ is the largest pseudometric on $X$ which satisfies (i)
and (ii).

\medskip

\noindent In fact, $k_X$ is the integrated form of the {\it
ROYDEN INFINITESIMAL\\ PSEUDOMETRIC:}
$$\forall x\in X,\,\,\forall v\in T_xX,$$
$$K_X(x,v) = \inf\{r^{-1}\,\vert\,
\exists \varphi\in {\rm HOL}\,(\Delta_r,\, X)\,:\,\,0 \longmapsto
x,\,\,\,\partial/\partial x \longmapsto v\}$$

\bigskip

\noindent {\bf REMARK}: Actually, any holomorphic map $\varphi: Y
\to X$ of complex spaces is a contraction: $\varphi^*(k_X)\le
k_Y$.

\bigskip

\noindent {\bf DEFINITION}: A complex space $X$ is called
$KOBAYASHI\,\,\,\,HYPERBOLIC$\\ if $k_X$ is a metric, i.e.
$k_X(p,\,q)=0 \Leftrightarrow p=q.$

\bigskip

\noindent {\bf EXAMPLES}: 1. The Kobayashi pseudometric on a
projective space $\PP^n$ (on an affine space $\C^n$, on a complex
torus $\C^n/\Lambda$, respectively) vanishes identically.

\bigskip

\noindent  2. Any bounded domain in $\C^n$ is hyperbolic. In
particular, the Teichm\"uller space $X=T_{g,n}$ is hyperbolic, and
the Kobayashi metric $k_X$ coincides with the Teichm\"uller
metric ({\bf Royden 1971}).

\newpage

\centerline{\underline{\bf Geometric function theory}}

\centerline{\underline{\bf on hyperbolic complex spaces}}

\vskip 0.3in

\noindent {\bf Schottky-Landau THEOREM:} {\it $\C\setminus
\{0,\,1\}$ is hyperbolic.}

\vskip 0.2in

\noindent Thus '$X$ is hyperbolic` means that the Schottky-Landau
Theorem holds for $X$.

\vskip 0.2in

\noindent {\bf Brody-Kiernan-Kobayashi-Kwack THEOREM:} {\it Let
$X$ be a compact complex space. Then the following are equivalent:

\medskip

\noindent (i) $X$ is Kobayashi hyperbolic;

\medskip

\noindent (ii) {\bf (Picard Theorem)} \newline $X$ does not
contain entire curves i.e. $\forall f\in {\rm HOL}\,(\C,\,
X),\,\,\,\,f= {\rm const};$

\medskip

\noindent (iii) {\bf (Big Picard Theorem)} \newline
$$\forall f\in {\rm HOL}\,(\Delta\setminus\{0\},\, X)\,\,\,\,\exists
{\bar f}\in {\rm HOL}\,(\Delta,\, X),\,\,\,\,
{\bar f}\vert (\Delta\setminus\{0\})=f;$$

\medskip

\noindent (iv) {\bf (Montel Theorem)} the space ${\rm
HOL}\,(\Delta,\,X)$ is compact.}

\vskip 0.1in

\noindent {\bf REMARK:} Actually, $X$ is hyperbolic
$\Longrightarrow $ for any complex space $Y$,\\ the space ${\rm
HOL}\,(Y,\,X)$ is compact.

\vskip 0.1in

\noindent {\bf STABILITY THEOREM:} {\it (a) {\bf (Brody
\cite{Br})} If a compact complex subspace $X\subset Z$ of a
complex space $Z$ is hyperbolic then there exists a hyperbolic
neighborhood $U\supset X$ in $Z$. Consequently, any complex
subspace $X'$ in $Z$ sufficiently close to $X$ is hyperbolic.

\smallskip

\noindent (b) {\bf (Zaidenberg \cite{Za2})} Suppose $Z$ is
compact and $X$, $X'$ are divisors in $Z$ such that $X$,
$Z\backslash X$ are hyperbolic. If $X'$ is sufficiently close to
$X$ then as well $X'$, $Z\backslash X'$ are hyperbolic.}

\vskip 0.1in

\noindent {\bf THEOREM:} {\it If $f\,:\,Y\to X$ is a covering
then $f$ is a local isometry: $f^*(K_X)=K_Y$. }

\bigskip

\noindent {\bf COROLLARY:}
{\it $X$ is hyperbolic
$\Leftrightarrow Y$ is hyperbolic.}

\vskip 0.1in

\noindent {\bf REMARK:} In particular, the Kobayashi metric on a
Riemann surface $X$ of non-exceptional type coincides with the
Poincar\'e metric of $X$.

\vskip 0.1in

\noindent {\bf COROLLARY (`The Covering Trick'):} {\it Let
$H\subset \PP^{n+1}$ be a hypersurface and $\pi\,:\,H\to \PP^n$ a
projection ramified over $X\subset \PP^n.$ If $H$ is hyperbolic
then $\PP^n\setminus X$ is hyperbolic as well.}

\newpage

\centerline{\underline{\bf 2. Kobayashi Problems}}

\vskip 0.2in

\noindent {\bf PROBLEM (Kobayashi 1970 \cite{Ko1}):} {\it Let $X$
be a (very) general hypersurface in $\PP^n$ of degree $d$. Is it
true that for any $d>>1,\,\,\,\,X$ resp., $\PP^n\setminus X$ is
hyperbolic? Is this true already for $d\ge 2n-1$ resp., $d\ge
2n+1$?}

\vskip 0.3in

\centerline{\underline{\bf Hyperbolicity of complements}}

\vskip 0.2in

\noindent The latter bound in the Kobayashi Problem is due to the
famous

\vskip 0.2in

\noindent {\bf Borel LEMMA (E. Borel 1898; Bloch 1926):} {\it
$\PP^n \setminus \{2n+1$ hyperplanes in general position$\}$ is
hyperbolic.}

\vskip 0.1in

\noindent In fact, this bound is optimal:

\vskip 0.1in

\noindent {\bf PROPOSITION (Zaidenberg \cite{Za1}):} {\it The
complement of a general hypersurface of degree $d\le 2n$ in
$\PP^n$ is not hyperbolic.}

\vskip 0.1in

\noindent {\bf THEOREM (Siu-Yeung \cite{SY2}, $d\ge 10^6$; El
Goul \cite{EG}):}  {\it For a very generic curve $C\subset \PP^2$
of degree $d\ge 13$ the complement $\PP^2\setminus C$ is
hyperbolic.}

\vskip 0.2in

\noindent Thus the hyperbolicity of the complement of a generic
plane curve of degree $d$ is still to be established for $5\le d
\le 12$. Starting with dimension 3, the Kobayashi problem on
hyperbolicity of complements is open. As for examples,  an
irreducible plane curve of degree 30 with hyperbolic complement
was constructed by
{\bf Azukawa-Suzuki 1980}. \\
Later on, examples in smaller degrees were found by {\bf
Grauert-Peternell, Carlson-Green, Adachi-Suzuki,
Dethloff-Schumacher-Wong, Dethloff-Zaidenberg, Berteloot-Duval,}
e.a. Examples of projective hypersurfaces in any dimension with
hyperbolic complements were suggested by {\bf Masuda-Noguchi
1996, Fujimoto 2000} e.a.

\vskip 0.1in

\noindent {\bf THEOREM:} {\it (a) {\bf (Zaidenberg \cite{Za2})}
For every $d\ge 5$ there exists an irreducible smooth curve $C$ in
$\PP^2$ of degree $d$ with hyperbolic complement. However, there
is no such curve of degree $d\le 4$.

\smallskip

\noindent (b) {\bf (Green \cite{Gr1}}; {\rm see also} \cite{Ba,
Za2, ErSo}{\bf )} Let $D=\bigcup_{i=1}^{2n+1} D_i$ be a union of
hyperbolic hypersurfaces in general position in $\PP^n$. Then for
any (smooth) hypersurface $D'\subset \PP^n$ close enough to $D$,
the complement $\PP^n\setminus D'$ is hyperbolic.

\smallskip

\noindent (c) {\bf (Shiffman-Zaidenberg \cite{ShZa1})} There exist
algebraic families of hyperbolic hypersurfaces of degree
$d=(2n-1)^2+2n$ in $\PP^n$ with hyperbolic complements.}

\bigskip

\noindent In particular (c) provides algebraic families of curves of degree
13 in $\PP^2$, of surfaces of degree 31 in $\PP^3$, and so forth,
whose complements are complete hyperbolic and hyperbolically
embedded in projective space. These examples are of Fermat-Waring
type (see below).

\newpage

\centerline{\underline{\bf Hyperbolicity
of projective hypersurfaces}}

\vskip 0.3in

\noindent {\bf THEOREM (McQuillan \cite{MQ}, $d\ge 36$;
Demailly-El Goul \cite{DEG}, $d\ge 21$):}\\
{\it A very general surface of degree $d \ge 21$ in $\PP^3$ is
Kobayashi hyperbolic.}

\vskip 0.3in

\noindent Thus the hyperbolicity of a generic surface in $\PP^3$
of degree $d$ is still to be established for $5\le d \le 20$.
Starting with dimension $n=4$ the Kobayashi problem on
hyperbolicity of a generic hypersurface in $\PP^n$ is
open\footnote{{\bf Y.-T. Siu} (unpublished) announced a positive
solution.}.

\vskip 0.3in

\noindent {\bf REMARK:} By the Stability Theorem, if a degree $d$
hypersurface $X$ in $\PP^n$ is hyperbolic then so is any
hypersurface $X'$ in $\PP^n$ of the same degree, sufficiently
close to $X$. Thus such hypersurfaces are parametrized by an open,
in the classical topology, subset $\cal H$ in $\PP^{N-1}$, where
$N= \binom{n+d}{d}$. The Kobayashi Problem asks, in particular,
whether $\cal H$ contains a Zariski open subset. What can be said
about the boundary of $\cal H$?

\vskip 0.5in

\centerline{\underline{\bf Methods employed:}}

\vskip 0.3in

\noindent 1. Value distribution theory;
the Nevanlinna-Cartan theorem\\
{\bf (Green 1975, Masuda-Noguchi 1996, McQuillan-Brunella 1998,
e.a.)}

\bigskip

\noindent 2. Meromorphic connections\\
{\bf (Siu 1987, Nadel 1989, Demailly-El Goul 1997, e.a.)}\\

\bigskip

\noindent 3. Jet differentials/jet metrics\\
{\bf (Bloch 1926, Green-Griffith 1979, Grauert 1989, Siu-Yeung 1996,
Demailly-El Goul 1997, Dethloff-Lu 1996, e.a.)} \\

\bigskip

\noindent 4. Algebraic (multi)foliations\\
{\bf (McQuillan 1997, Demailly-El Goul 1997, e.a.)}

\newpage

\centerline{\underline{\bf Examples}}

\medskip

\centerline{\underline{\bf of hyperbolic surfaces in $\PP^3$ of
degree $d$:}}

\vskip 0.2in

\noindent Below $a, a_i\in \C$ are generic.

\vskip 0.2in

\noindent {\bf Brody-Green 1977, $d=2k\ge 50$:}

$$X_0^{d}+X_1^{d}+X_2^{d}+X_3^{d}+a_1(X_0X_1)^k +
a_2(X_2X_3)^k=0\,.$$

\vskip 0.2in

\noindent {\bf Masuda-Noguchi 1996, $d=3e\ge 24$:}

$$X_0^{d}+X_1^{d}+X_2^{d}+X_3^{d}+a(X_0X_1X_2)^e=0\,,$$
$d=4e\ge 28$:
$$X_0^{d}+X_1^{d}+X_2^{d}+X_3^{d}+a(X_0X_1X_2X_3)^e=0\,,$$

\vskip 0.2in

\noindent and similar examples of hyperbolic hypersurfaces in
$\PP^n$ for all $n\ge 4$; they are of degree $d\ge 192$ for $n=4$
and $d\ge 784$ for $n=5.$

\vskip 0.2in

\noindent {\bf Khoai 1996, $d\ge 22$:}

$$X_0^{d}+X_1^{d}+X_2^{d}+X_3^{d}
+aX_0^{\alpha_0}X_1^{\alpha_1}X_2^{\alpha_2}=0$$ where
$\alpha_i\ge 7, \,\,\,\,\alpha_0+\alpha_1+\alpha_2=d\,.$

\vskip 0.2in

\noindent {\bf Nadel 1989, $d=6e+3\ge 21$:}

$$X_0^{6e}(X_0^3+aX_1^3)+X_1^{d}+X_2^{d}
+X_3^{6e}(X_3^3+aX_1^3)=0\,.$$

\vskip 0.2in

\noindent {\bf El Goul 1996, $d\ge 14$:}

$$X_0^{d-2}(X_0^2+aX_1^2)+X_1^{d}+X_2^{d}
+X_3^{d-2}(X_3^2+aX_1^2)=0\,.$$

\vskip 0.2in

\noindent {\bf Siu-Yeung 1996, Demailly-El Goul 1997, $d\ge 11$:}

$$X_0^d+X_1^{d}+X_2^{d}+X_3^{d-2}(a_0X_0^2+a_1X_1^2+a_2X_2^2+a_3X_3^2)=0\,.$$

\vskip 0.2in

\noindent {\bf J. Duval 2000 \cite{Du1}, Fujimoto
2000\footnote{Also {\bf Shirosaki 2000}, $d=2k\ge 10$ \cite{Shr1,
Shr2}.} \cite{Fu}, $d=2k\ge 8$:}
$$Q(X_0,X_1,X_2)^2-P(X_2,X_3)=0\,,$$
where $P,\,Q$ are general homogeneous forms of degree $d$
resp., $k$. The idea behind this example is quite simple. Let $S$
be a surface in $\PP^3$ as in the example.
A resolution $f: S'\to \PP^1$ of the meromorphic map
$(X_2:X_3) : S \dashrightarrow \PP^1$ has reducible general
fibers and so admits a Stein factorization $S'\stackrel{\varphi}
\to\Gamma\stackrel{2:1}\to\PP^1$, where $\Gamma$ is ramified over
the zeros of $P$. As $k\ge 4$ then $\Gamma$ and every fiber of $\varphi$
are irreducible curves of genera $\ge 2$.
Hence by Picard's Theorem any holomorphic map $\C\to S'$
is constant, and as well any holomorphic map $\C\to S$ is.
Thus by Brody's Theorem $X$ is hyperbolic.

\vskip 0.3in

\bigskip

\noindent Examples of hyperbolic hypersurfaces in any dimension were given
by {\bf Masuda-Noguchi 1996, Khoai 1996, Siu-Yeung 1997, Fujimoto
2000 and Shiffman-Zaidenberg 2000}. The best degree asymptotic
in these examples is achieved by the following
Fermat-Waring type hypersurfaces.

\bigskip

\noindent {\bf THEOREM (Siu-Yeung \cite{SY1}, $d=16(n-1)^2$,
Shiffman-Zaidenberg \cite{ShZa2}):}\\
{\it Let  $d\ge (m-1)^2$, $m\ge 2n-1$.  Then for generic linear
functions $h_1,\dots,h_m$ on $\C^{n+1}$, the hypersurface $X
\subseteq\PP^n$ with equation
$$\sum_{j=1}^{m} h^d_j=0$$
is hyperbolic. In particular there exist algebraic families of
hyperbolic hypersurfaces in $\PP^n$ of degree $d=4(n-1)^2$.}

\vskip 0.3in

\centerline{\underline{\bf 3. Symmetric powers of curves as
hyperbolic hypersurfaces in $\PP^3$ and $\PP^4$}}

\bigskip
\centerline{\underline{(after \cite{ShZa1}, \cite{CiZa})}}

\vskip 0.5in

\noindent {\bf DEFINITION:} A smooth projective curve $C$ is called
{\it HYPERELLIPTIC} resp., {\it BIELLIPTIC} if there exists a
$2:1$ morphism
$C \to \PP^1$ resp., $C\to E$ where $E$ is an elliptic curve.

\vskip 0.3in

\noindent {\bf LEMMA:} {\it (a)
The $n$-ths symmetric power $C^{(n)}$ of a generic
genus $g$ curve $C$
is hyperbolic if and only if $g\ge 2n-1$.

\noindent (b) The symmetric square $C^{(2)}$ of a
curve $C$ is hyperbolic iff $C$ is neither hyperelliptic nor
bielliptic. In particular, $C^{(2)}$ is hyperbolic for a genus
$g\ge 3$ curve $C$ with general moduli.}

\vskip 0.3in

\noindent {\bf REMARKS:} (a) A genus 2 curve is hyperelliptic.

\noindent (b) A non-hyperelliptic genus 3 curve $C$ is a smooth
plane quartic, and vice versa; it is bielliptic iff the group
AUT$(C)\subset $PSL$(3; \C)$ is of even order, that is, iff this
group contains an involution.

\vskip 0.3in

\noindent {\bf THEOREM \cite{ShZa1}:} {\it Let $C$ be a genus $g\ge 3$ curve
with general moduli, $C^{(2)}$ be its symmetric square embedded in
$\PP^5$, and $S$ be a general projection of $C^{(2)}$
in $\PP^3$. Then $S$ is a hyperbolic surface of degree $d\ge 16$.}

\newpage

\centerline{\underline{\bf Example of degree 16:}}

\vskip 0.3in

\noindent Let $C$ be a non-bielliptic
smooth plane quartic; e.g.
$$C=\{x^4-xz^3-y^3z=0\}\subset \PP^2.$$
Then the symmetric square $C^{(2)}$ of $C$ is a hyperbolic
surface in the symmetric square $S^{(2)}\PP^2$ of $\PP^2.$ In
turn, $S^{(2)}\PP^2$ can be embedded in $\PP^5$ as the cubic 4-fold
$$rst + 4uvw - r^2w - s^2u - t^2v =0$$
\noindent under $\varphi\,:\,S^{(2)}\PP^2\hookrightarrow \PP^5$,
$$\{(x:y:z),\,(x':y':z')\} \stackrel{\varphi}\longmapsto
(xy'+x'y:yz'+y'z:xz'+x'z:xx':yy':zz').$$

\medskip

\noindent This provides an embedding
$\varphi\,:\,C^{(2)}\hookrightarrow \PP^5$ of degree 16. A generic
projection to $\PP^3$ of the image $\varphi(C^{(2)})$  with
center at a generic line in $\PP^5$ is a hyperbolic surface in
$\PP^3$ of degree 16.

\bigskip

\noindent This surface is singular. By the Stability Theorem its
general small deformation is a smooth hyperbolic surface of degree
16 in $\PP^3$.

\vskip 0.2in

\noindent {\bf REMARK:} There are examples ({\bf
Kaliman-Zaidenberg} \cite{KaZa}) of non-hyperbolic singular
projective surfaces $X$ with a smooth hyperbolic normalization
$\tilde X \to X$.

\vskip 0.5in

\centerline{\underline{\bf Hyperbolic 3-folds in $\PP^4$
birational to the symmetric cube of a curve}}

\vskip 0.3in

\noindent {\bf THEOREM {\bf (Ciliberto-Zaidenberg \cite{CiZa})}:}
{\it Let $C$ be a  curve of genus $g\ge 7$ with general moduli,
and let $C^{(3)}$ be its symmetric cube embedded in $\PP^n,\,\,n\ge 7$.
If $T$ is a generic projection of $C^{(3)}$ in $\PP^4$ then $T$
is a hyperbolic 3-fold in $\PP^4$.
The same conclusion holds for certain special
embeddings $C^{(3)}\hookrightarrow\PP^7$, where $C$ is a general
plane quintic ($g=6$). The minimal degree of such a hyperbolic
3-fold in $\PP^4$ is 125.}
\\

\newpage


\centerline{\underline{\bf 4. Deformation method and more
hyperbolic surfaces in $\PP^3$}}

\bigskip

\centerline{\underline{(after \cite{ShZa3, ShZa4})}}

\vskip 0.3in

\noindent In Example 3 below we  give a simple
construction of hyperbolic surfaces
in $\PP^3$
of any given degree $d\ge 8$. We need the following notion.

\vskip 0.2in

\noindent {\bf DEFINITION:} Let $M$ be a complex manifold. A {\it
BRODY CURVE\/} in $M$ is an entire curve $\varphi:\C\to M$ such
that $\|\varphi'(\zeta)\|$ is bounded by $\|\varphi'(0)\|=1$,
where the norm is computed with respect to a hermitian metric on
$M$.

\bigskip

\noindent {\bf THEOREM (Brody \cite{Br})} {\it A compact complex
manifold $M$ is hyperbolic iff it does not contain Brody
curves.}

\vskip 0.3in

\noindent {\bf DEFORMATION METHOD\/:} Let $\{X_t\}=\langle X_0,
X_{\infty}\rangle$ be a linear pencil of degree $d$ surfaces in
$\PP^3$ generated by $X_0, X_{\infty}$, where $X_{\infty}$ is
general and $X_0$ is singular with a double curve
$\Gamma\subseteq $\,sing $X_0$ i.e., $X_0$ has two branches at
general points of $\Gamma$. Then sufficiently small deformations
$X_t$ of $X_0$ are hyperbolic provided that $X_0$, while not
hyperbolic, satisfies the following 'hyperbolic non-percolation`
property with respect to a certain divisor $D\supseteq \Gamma\cap
X_{\infty}$ on $\Gamma$, which will be precised later on.

\vskip 0.3in

\noindent {\bf DEFINITION:} Let $D$ be a divisor on $\Gamma$. We
say that $X_0\backslash \Gamma$ {\it has the property of
HYPERBOLIC NON-PERCOLATION through\/} $D$ if there is no Brody
curve $\varphi:\C\to (X_0\backslash \Gamma)\cup D$.

\vskip 0.3in

\noindent {\bf INDICATION FOR THE DEFORMATION METHOD:}\\
Assume that, for
a sequence $t_n\to 0$, the
surfaces $X_{t_n}$ are not hyperbolic. By Brody's Reparametrization Lemma,
there exists a sequence of Brody curves $f_{n}:\C\to
X_{t_n}$ converging to a Brody curve $f:\C\to X_0$. By Hurwitz's
Theorem then either
\begin{itemize}
\item $f(\C)\subset
\Gamma\backslash \{p_j\}$, where the $p_j$ are the points of
$\Gamma$ where $X$ has 3 or more local branches, or
\item $f(\C)\subset (X_0\backslash \Gamma)\cup D$ with
$D$ consisting of $\Gamma\cap X_{\infty}$ and the points of $\Gamma$
where $X_0$ is unibranched.
\end{itemize}
Indeed, suppose that $f(\C)$ passes through a point $p\in
\Gamma\backslash X_{\infty}$, and $X_0$ has $k$ branches at $p$.
Locally in a small neighborhood $U$ of $p$, $X_t$ can be given by
an equation of the form
$\varphi_1\cdot\ldots\cdot\varphi_k+t\varphi_{\infty}=0$, where
$\varphi_i(p)=0$ and $\varphi_{\infty}$ does not vanish in $U$.
Let $\Delta_{\varepsilon}$ be a small disc in $\C$ with center
$a$ such that $f_n(\Delta_{\varepsilon})\subseteq U$ for all
$n>>1$ and $f(a)=p$. Then $\varphi_i\circ f_n$ does not vanish in
$\Delta_{\varepsilon}$ while $\varphi_i\circ f$ does. By
Hurwitz's Theorem $\varphi_i\circ f\equiv 0$ ($i=1,\ldots,k$), as
needed.

\smallskip

\noindent Consequently, if \begin{itemize}
\item $\Gamma\backslash \{p_j\}$ is hyperbolic, and
\item $X_0\backslash \Gamma$ has
the property of hyperbolic non-percolation through $D$,
\end{itemize}
then all sufficiently small deformations $X_t$ of $X_0$ are hyperbolic.

\newpage

\noindent {\bf EXAMPLE 1 \cite{ShZa3}:} We let $X_0$ be the union
of 15 hyperplanes $\{l_j=0\},\,\, j=1,\dots,15,$ in general
position in $\PP^3$ and $X_{\infty}=3Q$, where $Q=\{q=0\}$ is a
generic quintic surface. We consider the linear pencil $\langle
X_0, X_{\infty}\rangle$ of surfaces $X_t$ in $\PP^3$, where
$$\textstyle X_t = \left\{\prod_{j=1}^{15}l_j + t
q^3=0\right\}\,.$$ Then the surface $X_t$ is hyperbolic for all
sufficiently small $t\ne 0$.

\vskip 0.3in

\noindent {\bf CONJECTURE:} {\it A general small deformation of 6
planes in general position in $\PP^3$ is a hyperbolic sextic
surface.}

\bigskip

\noindent Actually,
it is enough to show that the complement of the union $L$ of
general 5
lines in $\PP^2$ has hyperbolic non-percolation property with
respect to $D=L\cap C$, where $C$ is a general sextic curve.

\vskip 0.3in

\noindent {\bf EXAMPLE 2 \cite{ShZa3}:} There exists an octic
surface $X_0$ in $\PP^3$ with a double curve $\Gamma$ such that
the normalization of $X_0$ is smooth and is a simple abelian
surface $A$. Then for a general octic surface $X_{\infty}$ in
$\PP^3$, all sufficiently small deformations $X_t\in\langle
X_0,\,X_{\infty}\rangle$ ($t\neq 0$) of $X_0$ are hyperbolic.

\vskip 0.3in

\noindent {\bf EXAMPLE 3 \cite{ShZa4}:} We let $X_0$ be the union
of two general cones $CF_i=\langle a_i, F_i\rangle$ in $\PP^3$
with vertices $a_i$ ($a_1\neq a_2$) over generic plane curves
$F_i$ of degree $d_i\ge 4,\,\,i=1,2$. Then for a general surface
$X_{\infty}$ in $\PP^3$ of degree $d=d_1+d_2$, all sufficiently
small deformations $X_t\in\langle X_0,\,X_{\infty}\rangle$
($t\neq 0$) of $X_0$ are hyperbolic. In suitable coordinates
$(Z_0:\ldots :Z_3)$ in $\PP^3$, $X_t$ is given by the equation
$$f_1(Z_0,Z_1,Z_2)f_2(Z_1,Z_2,Z_3)+tf_{\infty}(Z_0,Z_1,Z_2,Z_3)=0\,,$$
where $F_i=\{f_i=0\}$, $i=1,2$, and $X_{\infty}=\{f_{\infty}=0\}$.

\vskip 0.3in

\noindent {\bf INDICATION:} We let $\Gamma=CF_1\cap CF_2$ be the
double curve of $X_0$. By our genericity assumption, the
projection $\pi_{a_i}:\Gamma\to F_i$ with center $a_i$ has degree
$d_j\ge 4$ ($j\neq i$) and only simple ramifications. Hence every
fiber of $\pi_{a_i}|\Gamma$ contains at least 3 points. If nearby
surfaces $X_t$ were not hyperbolic, then one would find a
sequence of Brody curves $f_{t_n}:\C\to X_{t_n}$ converging to a
Brody curve $f:\C\to X_0$. Let $f(\C)\subseteq CF_i$
($i\in\{1,2\}$). Since $F_i$ has genus $g_i\ge 3$, by Picard's
Theorem $\pi_{a_i}\circ f : \C \to F_i$ is constant. Thus
$f(\C)\subseteq l$, where $l\cong \PP^1$ is a projective line
through $a_i$ and a point $b\in F_i$.

Moreover, by Hurwitz's Theorem, $f(\C)$ does not meet $\Gamma$
except, maybe, at points of $X_{\infty}\cap \Gamma$. But a
general surface $X_{\infty}$ does not pass through the
ramification points of the projection $\pi_{a_i}: \Gamma\to F_i$
and meets any fiber of this projection in at most one point. It
follows that $\Gamma \backslash X_{\infty}$ contains at least 3
points of $l$. Hence by Picard's Theorem $f: \C\to l\backslash
(\Gamma \backslash X_{\infty})$ must be constant. This is a
contradiction.

\vskip 0.3in

\noindent {\bf EXAMPLE 4 (Duval 2004 \cite{Du2}) :} The same
deformation-nonpercolation methods, applied iteratively in a
clever way, provide an example of a degree 6 hyperbolic surface in
$\PP^3$.

\newpage

\centerline{\underline{\bf 5. Algebraic hyperbolicity}}

\vskip 0.2in

\noindent {\bf DEFINITION}: A projective variety is said to be
{\it ALGEBRAICALLY HYPERBOLIC} if it does not contain rational
curves or images of abelian varieties (in particular, it does not
contain elliptic curves).

\bigskip

\noindent {\bf REMARK:} Clearly, if $X$ is hyperbolic then it is
algebraically hyperbolic. As for the converse, no (projective)
counter-example is known. Moreover, basing on the {\bf Bloch
Conjecture} (see below), {\bf Green} and {\bf Griffiths} proposed
the following one.

\bigskip

\noindent {\bf Green-Griffiths CONJECTURE 1979 :}\\
{\it If $X$ is a projective surface of general type then the
image of any entire curve $f:\C\to X$ is contained in a closed
algebraic curve on $X$.}

\bigskip

\noindent Thus, $f(\C)$ has to be contained in a rational or
elliptic curve on $X$. Clearly, the number of these curves of any
given degree is finite. {\bf Bogomolov's CONJECTURE 1975} says
that the total number of such curves in $X$ is as well finite.
Some partial results were obtained by {\bf Bogomolov 1977, Lu-Yau
1990, Lu-Miyaoka 1997, McQuillan 1998, Brunella 1999 e.a.} For
instance {\bf (Grant 1986)}, the Green-Griffiths conjecture is
true for surfaces with irregularity $q(X)\ge 2$.

\bigskip

\noindent It is well known and elementary that on any
hypersurface of degree $d\le 2n-3$ in $\PP^n$ there are
projective lines. The number of these lines on a general
hypersurface of degree $2n-3$ is finite (e.g., there are $27$
lines on a smooth cubic surface in $\PP^3$). A smooth quartic
surface $S$ in $\PP^3$ (and more generally, any algebraic $K3$
surface) if does not contain lines, however contains rational
curves and a family of elliptic curves {\bf
(Mumford-Bogomolov-Mori-Mukai 1981)}. The rational (respectively,
elliptic) curves on $S$ arise e.g., as
 sections by tritangent (respectively, bitangent) planes of $S$.

\bigskip

\noindent {\bf THEOREM (Clemens \cite{Cl}, Ein \cite{Ei}, Xu \cite{Xu},
Voisin \cite{Vo}, Pacienza \cite{Pa}): }\\ {\it Let $X$ be a
very generic hypersurface of degree $d$ in $\PP^n$.

\medskip

\noindent (a) If $d\ge 2n-1$ then $X$ is algebraically
hyperbolic. In particular, a very generic surface in $\PP^3$ of
degree $d\ge 5$ is.

\medskip

\noindent (b) For $n\ge
4$ and $d\ge 2n-2$, $X$ does not contain rational curves,
and for $n\ge 6$
and $d= 2n-3$, the only
rational curves on $X$ are lines.}

\bigskip

\noindent {\bf REMARKS:} (a) By the {\bf Clemens CONJECTURE}, a
general quintic threefold $X\subset \PP^4$ ($d=2n-3$) contains
only a finite number of rational curves of any given degree. The
number of these curves is predicted by the Mirrow Symmetry.

\medskip

\noindent (b) If, for a given $(n,d)$, there exists a hyperbolic
degree $d$ hypersurface in $\PP^n$ then, clearly, a very generic
such hypersurface is algebraically hyperbolic. Thus examples of
hyperbolic hypersurfaces provide an alternative potential approach
to (a) in the above theorem. However, e.g., for $n=3, \,d=5$, so
far we are lacking such examples.

\newpage

\newpage

\centerline{\underline{\bf 6. Bloch Conjecture}}

\vskip 0.3in

\noindent {\bf Bloch 1926} formulated a conjecture (and an idea of
the proof) concerning hyperbolicity of subvarieties of abelian
varieties $A=\C^n/\Lambda$ and of the complements of divisors
$A\setminus D$. The simplest result in this direction is the
following one.

\vskip 0.3in

\noindent {\bf THEOREM (Green 1978 \cite{Gr2}) :} {\it A closed
subvariety $X$ of a compact complex torus $\C^n/\Lambda$ is
hyperbolic iff it does not contain shifted subtori.}

\vskip 0.2in

\noindent Indeed, according to Brody's Theorem, $\C^n/\Lambda$ is
not hyperbolic iff there exists a Brody curve $f:\C\to X$.
Clearly, the covering Brody curve $\tilde f:\C\to \C^n$ (w.r.t.
the Euclidean metric on $\C^n$) is an affine linear map, hence the
closure of the image $f(\C)$ in $X$ contains a shifted subtorus.

\vskip 0.3in

\noindent {\bf REMARK:} A generic complex torus (resp., a generic
Abelian variety, resp., a generic Jacobian variety) is simple.

\vskip 0.3in

\noindent {\bf COROLLARY:} {\it Let  $T=\C^n/\Lambda$ be a simple
complex torus, i.e. it does not contain any proper subtorus of a
positive dimension. Then any proper subvariety $V$ of $T$ is
hyperbolic. For instance, the theta-divisor $\Theta$ in a simple
abelian variety $T$ is hyperbolic.}

\vskip 0.3in

\noindent Many subsequent efforts were done to
fix the Bloch Conjecture in the full generality \\
{\bf (Lang 1966, Ax 1972, Ochiai 1977, Noguchi 1977, 1996,
Noguchi-Winkelmann 2003, Green 1978, Green-Griffith 1979,
Kawamata 1980, R. Kobayashi 1991, D. Abramovich 1994, McQuillan
1995, Siu-Yeung 1996, 2003, Demailly 1996, Dethloff-Lu 1996,
e.a.)}

\newpage

$\,$ \vskip 2in

\centerline{\bf A P P E N D I X}

\vskip 1in

\noindent  For reader's convenience we place below, as an
appendix, the survey article \cite{Za3} which was published in an
edition with a limited access.

\newpage

\begin{center}
{\bf Hyperbolicity in Projective Spaces\\
Mikhail Zaidenberg}\\
Universit\'e Grenoble I\\
Institut Fourier de Math\'ematiques\\
38402 St Martin d'H\`eres-cedex, France
\end{center}

\rightline {\sl To Professor Shoshichi Kobayashi} \rightline {\sl
on the occasion of his sixtieth birthday} \vskip 0.2in

In 1970 Sh. Kobayashi posed the following problems [Ko1]:

{\it Let $D$ be a generic hypersurface of degree $d$ in
${\PP}^n$, where $d$ is large enough with respect to $n$.}

{\bf I} {\it Is it true that $D$ is hyperbolic?}

{\bf II} {\it Is it true that the complement ${\PP}^n \setminus
D$ is hyperbolic and, moreover, hyperbolically embedded into
${\PP}^n$ ? Is this true for $d \ge 2n+1$ ?}

For $n=2$ (starting with $d=4$) the answer to {\bf I} is
classically known to be positive, while for $n \ge 3$ the problem
is open.

The answer to {\bf II} is unknown even for $n=2$. It is positive
for $n=1, d \ge 3$, and this is equivalent to the Montel Theorem.

Here we present a survey on the Kobayashi's Problems. Of course,
it does not pretend to be either exhaustive or original.

\begin{center}
\bf I The compact case
\end{center}

Let $\Pnd={\PP}^N$, where $N=\binom{n+d}{n}-1$, be the projective
space whose points parametrize (not necessarily reduced)
hypersurfaces of degree $d$ in ${\PP}^n$. Let ${\cal H}_{n,d}
\subset {\PP}_{n,d}$ be the subset corresponding to hyperbolic
hypersurfaces. To precise the meaning of "genericity" in {\bf I}
one could ask {\sl whether ${\cal H}_{n,d}$ contains a Zariski
open subset of ${\PP}_{n,d}$ for $d>>n$? Or, more generally,
whether the complement ${\PP}_{n,d} \setminus {\cal H}_{n,d}$ is
contained in a countable union of hypersurfaces in ${\PP}_{n,d}$
for $d>>n$ ?}

It is known that ${\cal H}_{n,d}$ is open (but probably empty) in
the classical Hausdorff topology of ${\PP}_{n,d}$ for any $n,d
\in {\N}$. This follows from the Brody's Stability Theorem [Br],
or, to be more precise, from the following version of it [Za1,4]:

\smallskip

{\bf Theorem I.1}  {\sl Let $M$ be a complex manifold and $X$ a
compact analytic subset of $M$. If $X$ is hyperbolic, then there
exists a neighborhood $U$ of $X$ in $M$, which is hyperbolically
embedded into $M$. Therefore, any compact analytic subset $X'$ in
$M$ close enough to $X$ is hyperbolic as well.

In particular, if $f: M \rightarrow S$ is a proper holomorphic
surjection onto a complex space $S$, then the subset of points in
$S$ that correspond to the hyperbolic fibers of $f$ is open.}

\smallskip

Let us give a sketch of proof.

Let $h$ be a fixed Hermitian metric on $M$. An entire curve $f:
{\C} \rightarrow M$ is called {\sl a Brody curve} iff $f$ is a
contraction with respect to the Euclidean metric in $\C$ and the
metric $h$ on M
 (i.e. $\mid df(z)\mid_h\le 1 \forall z\in{\C}$), and $\mid df(0)\mid_h=1$.

Let $\Delta_r$ be the open disc in $\C$ of radius $r$ centered at
the origin endowed with the metric $rh_r$, where $h_r$ is the
Poincar\'e metric in $\Delta_r$. It is easily seen that the
Euclidean metric in $\C$ is the limit of the metrics $rh_r$ as
$r\rightarrow \infty$. A holomorphic curve $f:
\Delta_r\rightarrow M$ is called {\sl a Brody curve} iff $f$ is a
contraction with respect to the metrics $rh_r$ in $\Delta_r$ and
$h$ in $M$, and $\mid df(0)\mid_h = 1$. By the Arzel\`a-Ascoli
Theorem any sequence $f_n: \Delta_n \rightarrow M$ of Brody
curves, whose images are contained in the same relatively compact
subset of $M$, has a subsequence converging to a Brody curve $f:
{\C}\rightarrow M$.

Let $\lbrace U_n\rbrace$ be a fundamental sequence of (relatively
compact) neighborhoods of the hyperbolic compact analytic subset
$X \subset M$. Suppose that there is no $n \in \N$ such that
$U_n$ is hyperbolically embedded into $M$. That means that the
inequality $K_{U_n}\ge ch$ for the Kobayashi-Royden pseudometric
$K_{U_n}$ on $U_n$ does not hold for any constant $c > 0$; in
particular, it does not hold for $c={1\over n}$. By the
definition of the Kobayashi-Royden pseudometric there exists a
sequence ${h_n: \Delta_n \rightarrow U_n}$ of holomorphic curves
such that $\mid dh_n(0) \mid > 1$. By the Brody Reparametrization
Lemma [Br] there exists a sequence of Brody curves $f_n:\Delta_n
\rightarrow
 U_n$, where $f_n (z)=h_n \circ \alpha_n (r_n z)$ for some $r_n <
1$ and $\alpha_n \in {\rm Aut}(\Delta_n)$. Passing to a convergent
subsequence, one obtains a limit Brody curve $f: {\C} \rightarrow
\bigcap U_n = X$, that contradicts our assumption on hyperbolicity
of $X$. $\bigcirc$
\\

So, the hyperbolicity of a hypersurface in $\PP^n$ is stable under
small deformations of the coefficients of the defining equation.
More generally, in any component of the Hilbert scheme of
projective varieties of a given degree and dimension, the set of
points which correspond to hyperbolic varieties is open in the
usual topology. We do not know {\sl when this set is non-empty;
whether, being non-empty, it must contain a Zariski open subset,
or at least an algebraic subvariety of small enough codimension.}

For $n=3$ R. Brody and M. Green [BrGre] gave examples of
one-parametric families of hyperbolic surfaces in ${\PP}^3$ of
any even degree $d=2k \ge 50$. Namely, the surfaces
$$D_{d,t}=\lbrace
{x_0}^{2k}+{x_1}^{2k}+{x_2}^{2k}+{x_3}^{2k}+t(x_0x_1)^k +
t(x_0x_2)^k = 0 \rbrace$$ (deformations of the Fermat surfaces
$F_{3,d}=D_{d,0}$) are hyperbolic for all but a finite number of
values of $t \in \C$. This means that for $d=2k \ge 50$ the set
${\cal H}_{3,d}$ is non-empty and contains a quasi-projective
rational curve $C=\lbrace D_{d,t}\rbrace$ (together with some
small classical neighborhood of it, as follows from the Stability
Theorem).

\smallskip

It is unknown whether for any $n \ge 4$ there exists a hyperbolic
hypersurface in $\PP^n$. J. Noguchi (private communication)
supposed that the Brody-Green construction should be available as
well in higher dimensions, and at least for $n=4$.

\smallskip

Notice that the Newton polyhedron of the Fermat hypersurface
$F_{n,d}$ of degree $d$ in $\PP^n$ is the standard simplex in
${\R}^{n+1}$; the monomials in the Fermat equation correspond to
its vertexes. Additional monomials in the Brody-Green example
correspond to the middle points of some edges of this simplex. So
the defining polynomials are fewnomials i.e., they contain few
monomials with respect to their degrees.

\smallskip

{\it Definition.} Let us say that a hypersurface $D=\lbrace p(x_0
,..., x_n)=0\rbrace$ of degree $d$ in $\PP^n$ is $k$-{\sl almost
simplicial} if any monomial of $p$ corresponds to a lattice point
in ${\R}^{n+1}$ with one of coordinates $\ge d - k$ (that is this
point is situated in a $k$-neighborhood of some vertex of the
$n$-simplex $\lbrace x_0 + ... + x_n = d\rbrace$ in $\R^{n+1}_+$).

\smallskip

The following statement is due to A. Nadel [Na]; its proof is
based on the Y.-T. Siu's version of the value distribution theory
for holomorphic curves in a complex manifold in presence of a
meromorphic connection.

\smallskip

{\bf Theorem I.2} {\sl For arbitrary $e \ge 3$ in the projective
space of all $k$-almost simplicial surfaces in ${\PP}^3$ of
degree $d=6e+3>4k+10$ there exists a quasiprojective subvariety
of dimension $4\binom{k+4}{4} - 1$, which consists of hyperbolic
smooth surfaces. In particular, ${\cal H}_{3,d}$ is non-empty for
any $d=6e+3\ge 21$.}

\smallskip

{\it Definition.} Let us say that a complex Hermitian manifold
$(X,h)$ is {\sl Brody hyperbolic} iff it does not contain any
Brody curve ${\C}\rightarrow X$, and {\sl Picard hyperbolic} iff
it does not contain any {non-constant} entire curve
${\C}\rightarrow X$.

\smallskip

The Picard Theorem says that ${\PP}^1 \setminus \lbrace 3\,\,
{\rm points} \rbrace$ is Picard hyperbolic. The Brody Theorem
[Br] states that for a compact manifold $X$ all three notions of
hyperbolicity (i.e., the Kobayashi hyperbolicity,  the Brody
hyperbolicity and  the Picard hyperbolicity) are equivalent.

M. Green [Gre4] noticed that a Brody curve ${\C}\rightarrow
{\T}^n$ in a complex torus ${\T}^n = {\C}^n/\Lambda$, where
$\Lambda$ is a maximal rank lattice in ${\C}^n = {\R}^{2n}$,
lifts to an isometric affine embedding ${\C}\rightarrow {\C}^n$.
Therefore, a closed subvariety $X\subset {\T}^n$ is (Brody)
hyperbolic iff it does not contain any shifted subtorus. The same
is true for any compact complex parallelizable manifold [HuWi].

More generally, Sh. Kobayashi [Ko2] established the following
fact.

\smallskip

{\bf Theorem I.3} {\sl Let $(X,h)$ be a Hermitian manifold with
non-positive holomorphic sectional curvature, and $f:
{\C}\rightarrow
 X$ be a Brody curve. Then $f$ is an isometric immersion, and its image
is totally geodesic.}

\smallskip

{\bf Problem I.1} {\sl Let the conditions of the above theorem be
fulfilled. Is it true that the closure $\overline{f({\C})}$ in
$X$ contains the image of a complex torus by a non-constant
holomorphic map, or at least any compact complex submanifold of
positive dimension?}

\smallskip

We notice that the rational curve ${\PP}^1$ and the simple complex
tori are the only known examples of compact complex manifolds with
totally degenerate Kobayashi pseudodistances that are minimal in
this class, i.e. that contain no closed subvarieties with this
property to be {\sl completely non-hyperbolic}. This motivates
the following

\smallskip

{\it Definition.} A compact complex space is said to be {\sl
algebraically hyperbolic} if it contains no image of a complex
torus by a non-constant holomorphic map.

\smallskip

In particular, an algebraically hyperbolic variety contains no
rational or elliptic curve. Clearly, a hyperbolic complex space
is also algebraically hyperbolic.

\smallskip

{\bf Problem I.2} {\sl Does algebraic hyperbolicity imply (Brody)
hyperbolicity, at least for projective varieties? In other words,
is it true that a compact complex space (a complex projective
variety) which possesses a Brody curve, must contain the image of
a complex torus under a non-constant holomorphic map?}

\smallskip

The following recent result of J.-P. Demailly and B. Shiffman
[DemSh] could be considered as an approximation to the positive
answer.

\smallskip

{\bf Theorem I.3} {\sl Let $X$ be a smooth projective variety,
$S$ a Stein manifold with $\dim S \le \dim X$, $f: S\rightarrow
X$ a holomorphic map,  $T$ a finite subset of $S$ and m a natural
number. Then there exists an exhaustive sequence $\Omega_1
\subset...\subset \Omega_k \subset...$ of Runge domains in $S$
and a sequence of holomorphic maps $f_k : \Omega_k \rightarrow
X_k$ such that, for any $k \in {\N}$, $\dim X_k = \dim S$ and at
each point $s \in T$ the $m$-jet of $f_k$ coincides with the
$m$-jet of $f$. If $S$ is an affine algebraic manifold, then
$f_k$ can be chosen to be regular.}

\smallskip

As a corollary, one gets the following algebraic definition of
the Kobayashi-Royden pseudometric $K_X$ on a projective variety
$X$: $$K_X (v) = \inf \lbrace K_{\tilde C} (v) \mid v\in
TC\rbrace ,$$ where $C$ runs over the set of all algebraic curves
in $X$ such that $v \in TX$ is a tangent vector to $C$, and
$K_{\tilde C}$ is the Poincar\'e metric of the normalization
$\tilde C$ of $C$. Furthermore, the Kobayashi pseudodistance $k_X
(x,y)$ on $X$ coincides with its algebraic analogue $d_X (x,y)$
suggested by J. Noguchi. Roughly speaking, the chains of
holomorphic discs in the definition of the Kobayashi
pseudodistance are replaced by chains of algebraic curves, and the
hyperbolic metrics of these curves are used instead of the
Poincar\'e metric in the disc.

\smallskip

An approach to Kobayashi's Problem {\bf I} is to divide it into
two parts: the above Problem I.2 on the equivalence of the Brody
hyperbolicity and the algebraic hyperbolicity for projective
varieties, as the first part, and as the second one the following

\smallskip

{\bf Problem I.3} {\sl Is it true that a generic projective
hypersurface of a large enough degree in $\PP^n$ is algebraically
hyperbolic?}

\smallskip

For $n=3$ the positive answer follows from the next recent result
of Geng Xu [Xu], which was conjectured by J. Harris and yields a
precision of an earlier one due to H. Clemens.

\smallskip

{\bf Theorem I.4} {\sl For any algebraic curve on a generic
surface $D \in {\PP}_{3,d}$ of degree $d \ge 5$ in ${\PP}^3$ the
following estimate holds: $$g(\tilde C) \ge {d(d - 3)\over 2} -2
\ge 3 , $$ where $g(\tilde C)$ is the genus of the normalization
$\tilde C$ of $C$. This bound is sharp, and for $d \ge 6$ the
curves of the minimal genus are sections of $D$ by tritangent
planes.

Therefore, for $d \ge 5$ a generic surface of degree $d$ in
${\PP}^3$ does not contain any rational or elliptic curve, and so
is algebraically hyperbolic.}

\smallskip

Observe that on a smooth quartic surface in ${\PP}^3$, and
moreover on any $K3$-surface, there exist a rational curve and a
linear pencil of elliptic curves (see [GreGri] and [MoMu]). Thus
such a surface is not algebraically hyperbolic. This shows that
the bound $d \ge 5$ above is sharp.

\smallskip

The proof of Theorem I.4 involves the Brill-Noether Theorem, and
thus the meaning of "genericity" in its formulation is more
extended than the genericity in Zariski sense. Namely, let ${\cal
AH}_{n,d} \subset {\Pnd}$ be the set of all algebraically
hyperbolic hypersurfaces. Then by Theorem I.4 for $d \ge 5$ the
complement ${\PP}_{3,d} \setminus {\cal AH}_{3,d}$ consists of a
countable number of proper algebraic subvarieties of
${\PP}_{3,d}$. There is no information about their mutual
position. In particular, the following problem seems to be
important.

\smallskip

{\bf Problem I.3} {\sl Is the locus ${\PP}_{3,d} \setminus {\cal
AH}_{3,d}$ closed in ${\PP}_{3,d}$ in the usual topology?}

\smallskip

Supposing this locus is not closed, there should exists a sequence
of non-algebraically hyperbolic surfaces $D_k$ in ${\PP}^3$
converging to an algebraically hyperbolic surface $D_0$. By the
stability of hyperbolicity, $D_0$ is not Brody hyperbolic;
indeed, otherwise for $k$ large enough $D_k$ would be hyperbolic
as well, and therefore algebraically hyperbolic. So, if the answer
to Problem I.3 were negative, then also the answer to Problem I.2
would be negative, and $D_0$ would be an example of an
algebraically hyperbolic surface which is not hyperbolic and hence
contains a Brody entire curve ${\C}\rightarrow D_0$.

\smallskip

A generic (in Zariski sense) hypersurface of degree $d \le 2n-3$
in $\PP^n$ contains a projective line (in particular, a smooth
cubic surface in ${\PP}^3$ contains exactly $27$ lines), thus is
not algebraically hyperbolic.

\smallskip

{\bf Question.} {\sl What is the maximal number $d=d(n)$ such
that $\Pnd \setminus {\cal AH}_{n,d}$ contains a Zariski open
subset of $\Pnd$?}

\smallskip

By the above remarks we have that $d(3) = 4$ and $d(n) \ge 2n-3$.

\smallskip

It is worthwhile also mentioning the following well known
problems:

\smallskip

{\sl Whether hyperbolicity (resp. algebraic hyperbolicity), or
even measure hyperbolicity of a compact complex manifold implies
that it is a projective variety of general type?}

\smallskip

The positive answer is known in the case of surfaces (see
[GreGri], [MoMu]).

\smallskip

A weaker property that could serve as a bridge between
hyperbolicity and algebraic hyperbolicity, is {\it algebraic
degeneracy}.

\smallskip

{\sl Definition.} One says that a complex space $X$ has the
property of {\it algebraic degeneracy} iff the image of any
non-constant entire curve ${\C} \rightarrow X$ lies in a proper
closed complex subspace of $X$. We mention {\it strong algebraic
degeneracy}, if this subspace is the same for all such curves.

\smallskip

Perhaps, it is worthwhile also to specify this notion by
restricting the class of curves under consideration to Brody
curves.

\smallskip

The Bloch Conjecture, proven by T. Ochiai, Y. Kawamata, and also
by M. Green and P. Griffiths, R. Kobayashi (see [RKo] for
references), states that {\sl an irregular projective variety $X$
(i.e. a variety with the irregularity $q(X) = h^{1,0}(X) > \dim
X$) has the property of algebraic degeneracy}. The above
restriction was weakened in the case of surfaces of general type
to $q(X) \ge 2$ by C. Grant [Gra1] (see also [Gra2], [HuWi], [Lu]
and St. Lu's report in this volume for some related results).

\smallskip

Another property, close to algebraic hyperbolicity, is finiteness
of the number of non-hyperbolic (resp. non-algebraically
hyperbolic) proper subvarieties. In the surface case this is
finiteness of the number of rational and elliptic curves, that
was proved by F. Bogomolov [Bo] for projective surfaces of
general type with ${c_1}^2 > c_2$ (see also [Lu]). H. Clemens
conjectured that the number of rational curves of any given
degree $d$ on a generic quintic threefold in ${\PP}^4$ is finite,
that was verified by N. Katz for $d \le 7$ (see [Xu]).

\vskip 0.1in

\begin{center}
\bf II The non-compact case
\end{center}

\smallskip

Denote by $\HE$ the subset of $\Pnd$ consisting of all
hypersurfaces of degree $d$ in $\PP^n$ with hyperbolically
embedded complements. Then $\HE$ is non-empty for any $d \ge
2n+1$; indeed, it contains the union $C_{n,d}$ of $d$ hyperplanes
in general position. This fact (modulo Kiernan's criterion of
hyperbolic embedding [Ki2]) goes back to E. Borel, A. Bloch, A.
Cartan and J. Dufresnoy (see [KiKo] for references). It was
reproved many times, for instance by M. Green [Gre2], E. Babets
[Ba] and others.

\smallskip

The bound $d \ge 2n+1$ for $\HE$ being non-empty should be sharp.
It is sharp for $n=2$; indeed, M. Green remarked in [Gre3] that
for any quartic curve $C$ in ${\PP}^2$ there exists a projective
line $l$ that intersects $C$ not more than in two points (an
inflectional tangent to $C$, a bitangent, a tangent in a singular
point, or a line passing through two singular points of $C$).
Thus ${\PP}^2 \setminus C$ is not hyperbolic; indeed, it contains
$l \setminus C \supset {\PP}^1 \setminus \lbrace 2\,\,{\rm
points} \rbrace$, and so the Kobayashi pseudodistance $k_{{\PP}^2
\setminus C}$ is degenerate along $l\setminus C$.

\smallskip

We do not know whether  for $d \le 2n$, $\HE$ is empty, although
we know [Za3] that its complement $\Pnd \setminus {\cal HE}_{n,d}$
contains a Zariski open subset.

\smallskip

{\bf Proposition II.1} {\sl For a generic (in Zariski sense)
hypersurface $D$ of degree $d \le 2n$ in $\PP^n$ and for any $k, 0
\le k \le d$, there exists a projective line $l$ that intersects
$D$ only in two points with multiplicities $k$ and $d - k$,
respectively. Thus, the pseudodistance $k_{{\PP^n} \setminus D}$
is degenerate along $l\setminus D$. For $d=2n$ the number of such
lines is finite.}

\smallskip

In contrast with $\HH\subseteq \Pnd$, the subset $\HE$ is never
open in the usual topology of $\Pnd$. For instance, for any $d \ge
2n+1$ the totally reducible hypersurfaces $C_{n,d} \in \HE$
considered above belong to the boundary of $\HE$. This follows
from the next simple observation [Za4]:

\smallskip

{\bf Proposition II.2} {\sl Any hypersurface $D_0$ in $\PP^n$ that
contains a projective line $l$, is the limit of a sequence of
hypersurfaces $\lbrace D_k \rbrace$ such that $l\cap D_k$
consists of a single point. Thus $\PP^n \setminus D_k$ is not
hyperbolic, and so $D_0 \in \overline{\Pnd \setminus \HE}$}.

\smallskip

However, in [Za4] a stability principle is obtained which can be
applied to fix the Kobayashi Problem II. Its proof follows the
line of the proof of Theorem I.1. It gives e.g., the following
result.

\smallskip

{\bf Theorem II.1} {\sl Let $M$ be a compact complex manifold and
$D$ a hypersurface in $M$. If $D$ and $M \setminus D$ are both
Brody hyperbolic, then $M \setminus D$ is hyperbolically embedded
in $M$. Moreover, all the above properties are preserved by small
deformations of the pair $(M,D)$.}

\smallskip

{\bf Corollary} {\sl $\HE \cap \HH$ is an open (but possibly
empty) subset of $\Pnd$ in the usual Hausdorff topology.}

\smallskip

Presumably,  for $d>>n$, the intersection $\HE \cap \HH$ contains
a Zariski open subset of $\Pnd$. This would imply the positive
answer to the both of the Kobayashi Problems.

\smallskip

To construct examples of hypersurfaces in $\HE \cap \HH$, one can
use the following generalization of the
Borel-Bloch-Cartan-Dufresnoy Theorem. It can be deduced from a
result of M. Green [Gre2], and it was proven by E. Babets [Ba] by
a different method.

\smallskip

{\bf Theorem II.2} {\sl The complement of the union of $2n+1$
smooth hypersurfaces in $\PP^n$ in general position is
hyperbolically embedded into $\PP^n$.}

\smallskip

In fact, this is true for any union of $2n+1$ hypersurfaces such
that the intersection of any $n+1$ of them is empty (A. Eremenko
and M. Sodin [ErSo]; a simplified proof has been recently done by
Min Ru). Using this theorem and Theorem II.1, one can easily
obtain the following

\smallskip

{\bf Corollary} {\sl If ${\cal H}_{n,k}$ is non-empty then  for
any $d \ge (2n+1)k$, the set $\HE \cap \HH$ is non-empty and
open.}

\smallskip

Indeed, by Theorem II.2 the union of any $2n+1$ smooth hyperbolic
surfaces in general position belongs to $\HE \cap \HH$.

\smallskip

In particular, from the existence of a hyperbolic surface in
${\PP}^3$ of degree $21$ [Na] it follows that ${\sl HE}_{3,d}
\cap {\sl H}_{3,d}$ is non-empty for any $d \ge 147 = 7\cdot 21$.

\smallskip

For $n=2$ a more refined version of the Stability Principle,
which uses {\sl absorbing stratifications} [Za4], leads to the
following result.

\smallskip

{\bf Theorem II.3} {\sl For any $d \ge 5$ the open set ${\sl
HE}_{2,d} \cap {\sl H}_{2,d}$ is non-empty i.e., there exists a
classically open set of smooth curves in ${\PP}^2$ of degree $d$
with hyperbolically embedded complements.}

\smallskip

The bound $d \ge 5$ is sharp, as follows from the remark of M.
Green mentioned above.

\smallskip

The first examples of smooth curves of any even degree $d \ge 30$
in ${\cal HE}_{2,d}$ were constructed by K. Azukawa and M. Suzuki
[AzSu] by the Brody-Green method [BrGre]. Note that, if
$B\subseteq {\PP}^2$ is the branch curve of a regular projection
to ${\PP}^2$ of a hyperbolic projective surface, then the
complement ${\PP}^2 \setminus B$ is a base of a hyperbolic
covering and so is hyperbolic. But the class of such curves is
rather restricted, as has been observed by F. Bogomolov, B.
Moishezon and M. Teicher. For instance, the number of cusps of the
branch curve of a generic projection to ${\PP}^2$ of a smooth
projective surface is divisible by $3$.

\smallskip

In a series of papers by M. Green, J. Carlson and M. Green, H.
Grauert and U. Peternell (see [Za4] for references) certain
sufficient conditions were found that ensure, for an irreducible
plane curve $C$ of genus $\ge 2$, the existence of a complete
Hermitian metric in the complement ${\PP}^2 \setminus C$ with
holomorphic sectional curvature bounded from above by a negative
constant. By Ahlfors Lemma this implies that ${\PP}^2 \setminus
C$ is hyperbolically embedded in ${\PP}^2$. Any curve satisfying
these conditions is singular and has degree $\ge 6$; the only
known examples are the dual curves to generic smooth plane curves
of degrees $d \ge 4$.

\smallskip

By Green-Babets Theorem II.2, the complement to a union of 5
smooth curves in ${\PP}^2$ in general position is hyperbolically
embedded in $\PP^2$. Therefore for $d \ge 5$ the set ${\cal
HE}_{2,d}$ contains a quasiprojective variety of positive
dimension. For instance, the quasiprojective subvariety $M=\lbrace
C_{2,5}\rbrace = \lbrace$the unions of $5$ lines in general
position in ${\PP}^2\rbrace$ of dimension $10$ is contained in
${\cal HE}_{2,5} \subset {\PP}_{2,5} = {\PP}^{20}$. Recently G.
Dethloff, G. Schumacher and P.-M. Wong [DetSchWo] have shown that
the complement to a union $C$ of $4$ plane curves in general
position is hyperbolically embedded in ${\PP}^2$ provided that
$\deg C\ge 5$ (see P.-M. Wong's report in this volume). This fact
can also be obtained by using a result of M. Green [Gre2], or the
one by Y. Adachi and M. Suzuki [AdSu1] (see Theorem II.6 below).

\smallskip

 Another result of [DetSchWo], conjectured by H. Grauert [Grau]
and obtained through the value distribution theory, is the
following

\smallskip

{\bf Theorem II.4} {\sl a) In the space of all unions of three
quadrics in ${\PP}^2$, there is an (explicitly defined) Zariski
open subset contained in ${\cal HE}_{2,6}$.

b) In the space of all unions of a line and two quadrics in
${\PP}^2$, there is a quasiprojective subvariety of codimension
$1$ contained in ${\cal HE}_{2,5}$. }

\smallskip

Let us mention a related criterion of hyperbolic embedding for the
complements of curves [Za2].

\smallskip

{\bf Proposition II.3} {\sl Let $C$ be a closed curve in a smooth
compact complex surface $M$. The complement $M \setminus C$ is
hyperbolically embedded in $M$ if and only if the curve $C
\setminus {\rm Sing}(C)$ is hyperbolic and the complement $M
\setminus C$ is Brody hyperbolic.}

\smallskip

The property of algebraic degeneracy of the complements of curves
was treated by T. Nishino and M. Suzuki [NiSu], Y. Adachi and M.
Suzuki [AdSu1,2]. The following results are worth mentioning.

\smallskip

{\bf Theorem II.5} ([NiSu]) {\sl Let $M$ and $C$ be as above. If
the logarithmic Kodaira dimension $\overline{k}(M \setminus C) =
2$, then any proper holomorphic map $f: {\C} \rightarrow  M
\setminus C$ is algebraically degenerate i.e., its image $f(C)$ is
contained in a closed curve $E$ in $M$.}

\smallskip

{\bf Theorem II.6} ([AdSu1]) {\sl If a reducible curve $C$ in
${\PP}^2$ consists of at least $4$ irreducible components which
do not belong to the same linear pencil, then there exists a curve
$A$ in ${\PP}^2$ such that the image of any non-constant entire
curve ${\C} \rightarrow {\PP}^2 \setminus C$ is contained in $A$.
Thus, ${\PP}^2 \setminus C$ has the property of strong algebraic
degeneracy.}

\smallskip

All possible exceptions here have been completely classified. For
some examples of degeneracy loci in the complements of irreducible
quartic curves see [Gre3]; see also [DetShuWo] for the reducible
case.

\smallskip

Another degeneracy principle has been used in the Babets' proof of
Theorem II.2 [Ba]. It states that, {\sl if $M$ is a compact
complex manifold and $D$ is a normal crossings divisor in $M$,
then, for a suitable complete Hermitian metric on $M\setminus D$,
every holomorphic differential in $M \setminus D$ with
logarithmic poles along $D$ is constant on any Brody curve ${\C}
\rightarrow M \setminus D$}. See also [Na] for an algebraic
degeneracy principle in presence of an ample (in Siu's sense)
meromorphic connection.

\smallskip

With evident changes, the notion of algebraic hyperbolicity can
be equally applied to affine or quasiprojective algebraic
varieties. This allows again to divide Problem II into two parts,
likewise Problem {\bf I} was divided above into Problems I.2 and
I.3.

\smallskip

{\bf Problem II.1} {\sl  Let $D$ be a hyperbolic hypersurface in
$\PP^n$ such that there exists a Brody curve ${\C} \rightarrow
\PP^n \setminus D$. Is it true that there exists a rational
projective curve $C$ in $\PP^n$ which has not more than two places
on $D$?}

\smallskip

{\bf Problem II.2} {\sl  Let ${\cal L}_{n,d} \subset \Pnd$ be the
locus of all hypersurfaces $D$ of degree $d$ in $\PP^n$ which
admit a rational curve $C$ as above. Is it true that,  for $d>>n$,
the complement $\Pnd \setminus {\cal L}_{n,d}$ contains a Zariski
open subset of $\Pnd$? Is the locus ${\cal L}_{n,d}$ closed in
$\Pnd$ with its Hausdorff topology?}

\smallskip

Next we pass to hyperbolicity properties of the complements to
hyperplanes in $\PP^n$. For hyperplanes in general position, the
following result is due to H. Fujimoto [Fu], M. Green [Gr], P.
Kiernan and Sh. Kobayashi [KiKo].

\smallskip

{\bf Theorem II.7} {\sl Let $D$ be a union of $n+k$ hyperplanes in
general position in $\PP^n$, where $k>0$. Then the image of any
non-constant entire curve ${\C} \rightarrow  {\PP^n} \setminus D$
is contained in a linear subspace of dimension $\le [{n\over
k}]$. This bound is sharp. Moreover the degeneracy locus is
contained in a finite union of the 'diagonal linear subspaces` of
dimension $n-k+1$ defined by $D$ in a canonical way. Thus $\PP^n
\setminus D$ has the property of strong algebraic degeneracy.}

\smallskip

For $k=2$ this gives the upper bound $[{n \over 2}]$ for the
dimension of the degeneracy locus. Observe that from the Borel
Lemma it  just follows the linear degeneracy, which means that any
non-constant entire curve in the complement to $n+2$ hyperplanes
in $\PP^n$ in general position is contained in a hyperplane. In
fact, the latter remains true without the assumption of general
position [Gre1]. For $k=n+1$, Theorem II.7 once again leads to the
Borel-Bloch-Cartan-Dufresnoy Theorem.

\smallskip

The bound $d \ge 2n+1$ for the hyperbolicity of $\PP^n \setminus
D$ is sharp, as is shown by the following result of V.E.
Snurnitsyn [Sn], which confirms a conjecture of P. Kiernan [Ki1].

\smallskip

{\bf Theorem II.8} {\sl For any union $D$ of $2n$ hyperplanes in
$\PP^n$ there exists a projective line which meets $D$ at most in
two points. Therefore, $\PP^n \setminus D$ is not hyperbolic.}

\smallskip

Some examples of unions of hyperplanes in non-general position
with hyperbolically embedded complements were given by P. Kiernan
[Ki1]. In [Za2] the following conditions for a finite union $D$
of hyperplanes in $\PP^n$ were considered:

\smallskip

(a) {\sl There does not exist a pair of points $x, y$ in $\PP^n$
such that each hyperplane in $D$ passes through at least one of
these points. In other words, there does not exist a projective
line $l=(x,y)$ which intersects the union of all those
hyperplanes in $D$ that do not contain $l$, in at most two
points.}

\smallskip

 (b) {\sl There does not exist a pair of points $(x, y)$
in $\PP^n$ such that each hyperplane in $D$ passes through
exactly one of these points. In other words, there does not exist
a projective line $l= (x,y)$ that intersects $D$ in at most two
points.}

\smallskip

If condition (b) fails then, clearly, the Kobayashi pseudodistance
$k_{\PP^n \setminus D}$ is degenerate along $l$. If (a) is
violated then the limit of $k_{\PP^n \setminus D}$ is degenerate
along $l$. The following criteria were obtained in [Za2, Sect.3].

\smallskip

{\bf Theorem II.9} {\sl Let $D$ be as above. The complement $\PP^n
\setminus D$ is hyperbolically embedded in $\PP^n$ if and only if
(a) holds. It is Picard hyperbolic if and only if (b) is
fulfilled. Furthermore, for $n=2$ (b) is equivalent to the
hyperbolicity of ${\PP}^2 \setminus D$.}

\smallskip

The latter had been conjectured by S. Iitaka.

\smallskip

Another criterion for the Picard hyperbolicity of complements of
hyperplanes has been recently obtained by Min Ru [Ru].

\smallskip

{\bf Theorem II.10} {\sl The complement $\PP^n \setminus D$ of a
finite union $D$ of hyperplanes in $\PP^n$ is Picard hyperbolic
if and only if, for any linear subspace $V$ in $\PP^n$ which is
not contained in $D$, the intersection $V \cap D$ contains at
least three distinct hyperplanes of $V$ that are linearly
dependent.}

\smallskip

The latter condition is obviously equivalent to (b). An algorithm
that allows to check it is given in [Ru]. To verify (b) one can
equally apply an algorithm similar to the simplex method, which
consists in passing from one pair of isolated intersection points
of $n$ hyperplanes in $D$ (if there is any such pair) to another
one.

\smallskip

In conclusion, let us mention the Lang Conjecture on equivalence
of Picard hyperbolicity and mordelleness (see [La]). For the
complements of hyperplanes, it was proven by P.-M. Wong and M. Ru
[WoRu] under the assumption of general position, and by M. Ru [Ru]
without this assumption.

\newpage

\begin{center}
{\bf References}
\end{center}

\smallskip

[AdSu1] Y. Adachi, M. Suzuki. {\sl On the family of holomorphic
mappings into projective space with lacunary hypersurfaces}. J.
Math. Kyoto Univ. 30 (1990), 451-458

[AdSu2] Y. Adachi, M. Suzuki. {\sl Degeneracy points of the
Kobayashi pseudodistances on complex manifolds}. Proc. Symp. Pure
Math. 52 (1991), P. 2, 41-51

[AzSu] K. Azukawa, M. Suzuki. {\sl Some examples of algebraic
degeneracy and hyperbolic manifolds}. Rocky Mountain J. Math. 10
(1980), 655-659

[Ba] V. A. Babets. {\sl Picard-type theorems for holomorphic
mappings}. Siberian Math. J. 25 (1984), 195-200

[Bo] F. A. Bogomolov. {\sl Families of curves on a surface of
general type}. Soviet Math. Dokl.18 (1977), 1294-1297

[Br] R. Brody. {\sl Compact manifolds and hyperbolicity}. Trans.
Amer. Math. Soc. 235 (1978), 213-219

[BrGre] R. Brody, M. Green. {\sl A family of smooth hyperbolic
surfaces in $\PP^3$}. Duke Math. J. 44 (1977), 873-874

[Co] M. Cowen. {\sl The method of negative curvature: the
Kobayashi metric on ${\PP}^2$ minus four lines}. Trans. Amer.
Math. Soc. 319 (1990), 729-745

[DemSh] J.-P. Demailly, B. Shiffman. {\sl Algebraic
approximations of analytic maps from Stein domains to projective
manifolds}, preprint (1992)

[DetSchWo] G. Dethloff, G. Schumacher, P.-M. Wong. {\sl
Hyperbolicity of the complements of plane algebraic curves},
preprint Math. G\"{o}tting. 31 (1992), 1-38

[ErSo] A. E. Eremenko, M. L. Sodin. {\sl The value distribution
for meromorphic functions and meromorphic curves from the point
of view of potential theory}. St. Petersburg Math. J. 3 (1992),
No. 1, 109-136

[Fu] H. Fujimoto. {\sl Families of holomorphic maps into
projective space omitting some hyperplanes}. J. Math. Soc. Japan
25 (1973), 235-249

[Gra1] C. Grant. {\sl Entire holomorphic curves in surfaces}.
Duke Math. J. 53 (1986), 345-358

[Gra2] C. Grant. {\sl Hyperbolicity of surfaces modulo rational
and elliptic curves}. Pacific J. Math. 139 (1989), 241-249

[Grau] H. Grauert. {\sl Jetmetriken und hyperbolische Geometrie}.
Math. Z. 200 (1989), 149-168

[Gre1] M. Green. {\sl Holomorphic maps into $\PP^n$ omitting
hyperplanes}. Trans. Amer. Math. Soc. 169 (1972), 89-103

[Gre2] M. Green. {\sl Some Picard theorems for holomorphic maps to
algebraic varieties}. Amer. J. Math. 97 (1975), 43-75

[Gre3] M. Green. {\sl Some examples and counterexamples in value
distribution theory}. Compos. Math. 30 (1975), 317-322

[Gre4] M. Green. {\sl Holomorphic maps to complex tori}. Amer. J.
Math. 100 (1978), 109-113

[GreGri] M. Green, P. Griffiths. {\sl Two applications of
Algebraic Geometry to entire holomorphic mappings}. In: "The
Chern Symposium 1979", Springer, N.Y. e.a. (1980), 41-74

[HuWi] A. Huckleberry, J. Winkelmann. {\sl Subvarieties of
parallelizable manifolds}, preprint (1992); to appear in Math.
Annalen

[Ki1] P. Kiernan. {\sl Hyperbolic submanifolds of complex
projective space}. Proc. Amer. Math. Soc. 22 (1969), 603-606

[Ki2] P. Kiernan. {\sl Hyperbolically imbedded spaces and the big
Picard theorem}. Math. Ann. 204 (1973), 203-209

[KiKo] P. Kiernan, Sh. Kobayashi. {\sl Holomorphic mappings into
projective space with lacunary hyperplanes}. Nagoya Math. J. 50
(1973), 199-216

[RKo] R. Kobayashi. {\sl Holomorphic curves into algebraic
subvarieties of an abelian variety}. Internat. J. Math. 2 (1991),
711-724

[Ko1] Sh. Kobayashi. {\sl Hyperbolic manifolds and holomorphic
mappings}. Marcel Dekker, 1970

[Ko2] Sh. Kobayashi. {\sl Complex manifolds with nonpositive
holomorphic sectional curvature and hyperbolicity}. Tohoku Math.
J. 30 (1978), 487-489

[La] S. Lang. {\sl Hyperbolic and Diophantine analysis}. Bull.
Amer. Math. Soc. 14 (1986), 159-205

[Lu] St.Sh.-Y.Lu. {\sl On meromorphic maps between algebraic
varieties with log-general targets}. Thesis, Harvard Univ.,
Cambridge, Mass., 1990

[MoMu] S. Mori, S. Mukai. {\sl The uniruledness of the moduli
space of curves of genus $11$}. In: "Algebraic Geometry
Conference (Tokyo-Kyoto 1982)", Lect. Notes in Math. 1016, 334-353

[Na] A. Nadel. {\sl Hyperbolic surfaces in ${\PP}^3$}. Duke Math.
J. 58 (1989), 749-771

[NiSu] T. Nishino, M. Suzuki. {\sl Sur les singularit\'es
essentielles et isol\'ees des applications holomorphes \`a
valeurs dans une surface complexe}. Publ. RIMS 16 (1980), 461-497

[Ru] M. Ru. {\sl Geometric and arithmetic aspects of ${\PP}^n$
minus hyperplanes}, preprint (1992)

[RuWo] M. Ru, P.-M. Wong. {\sl Integral points of $\PP \setminus
\lbrace 2n+1$ hyperplanes in general position$\rbrace$}. Invent.
Math. 106 (1991), 195-216

[Sn] V. E. Snurnitsyn. {\sl The complement of $2n$ hyperplanes in
${\C\PP}^n$ is not hyperbolic}. Matem. Zametki 40 (1986), 455-459
(in Russian)

[Xu] G. Xu. {\sl Subvarieties of general hypersurfaces in
projective space}, preprint (1992)

[Za1] M. Zaidenberg. {\sl The Picard theorems and hyperbolicity}.
Siberian Math. J. 24 (1983), 858-867

[Za2] M. Zaidenberg. {\sl On hyperbolic embedding of complements
of divisors and the limiting behavior of the Kobayashi-Royden
metric}. Math. USSR Sbornik 55 (1986), 55-70

[Za3] M. Zaidenberg. {\sl The complement of a generic
hypersurface of degree $2n$ in ${\C\PP}^n$ is not hyperbolic}.
Siberian Math. J. 28 (1987), 425-432

[Za4] M. Zaidenberg. {\sl Stability of hyperbolic imbeddedness and
construction of examples}. Math. USSR Sbornik 63 (1989), 351-361
\end{document}